\documentclass[11pt]{article}
\usepackage[T1]{fontenc}
\usepackage[margin=1in]{geometry}
\usepackage{amsmath,amsthm}
\usepackage{newtxtext}
\usepackage[subscriptcorrection]{newtxmath}
\usepackage{graphicx}
\usepackage[authoryear,round]{natbib}
\usepackage{microtype}
\usepackage{authblk}
\usepackage[hidelinks]{hyperref}
\graphicspath{{art/}}
\newtheorem{theorem}{Theorem}
\newtheorem{proposition}{Proposition}
\DeclareMathOperator{\E}{E}
\DeclareMathOperator{\pr}{pr}
\DeclareMathOperator{\var}{var}
\DeclareMathOperator{\supp}{supp}
\DeclareMathOperator{\law}{law}
\newcommand{\dd}{\,\mathrm{d}}
\newcommand{\dto}{\xrightarrow{d}}
\newcommand{\prto}{\xrightarrow{\mathrm{pr}}}
\newcommand{\ind}[1]{I(#1)}
\title{Nonsense nonlinear correlations between independent random walks}
\author{Liping Zhu\thanks{Email: \texttt{zhu.liping@ruc.edu.cn}.}}
\author{Zhongming Xie\thanks{Email: \texttt{xzm1@ruc.edu.cn}.}}
\affil{Institute of Statistics and Big Data, Renmin University of China\\
59 Zhongguancun Street, Haidian District, Beijing, China}

\date{}
\hypersetup{pdftitle={Nonsense nonlinear correlations between independent random walks},
 pdfauthor={Liping Zhu and Zhongming Xie}}
\begin{document}
\maketitle

\begin{abstract}
Sample distance correlation can remain large between independent random walks. For centred walks with independent, identically distributed increments of finite, positive variance, we prove convergence in distribution to the distance correlation of the joint occupation measure of two independent Brownian paths. This universal limit lies strictly between zero and one almost surely and has full support on the unit interval. We calculate the expected squared distance covariance of the limiting measure and show that ordinary permutation calibration rejects independence with probability tending to one at every fixed nominal level. Gaussian simulations illustrate the persistent dispersion of random-walk correlations and the failure of permutation calibration.
\end{abstract}


\section{Introduction}
\label{sec:introduction}

\citet{Yule1926} showed that independent time series can exhibit large sample Pearson correlations, a phenomenon he called ``nonsense correlation''. His random-walk experiments revealed much greater dispersion than samples of independent observations would suggest. \citet{Phillips1986} showed that these spurious correlations do not disappear as the sample size increases. \citet{Ernst2017} later quantified this phenomenon by deriving an analytic expression for the second moment of the empirical Pearson correlation between two independent Wiener processes. Numerical evaluation gives 0.240522, corresponding to a standard deviation of about 0.4904 and explaining the frequent occurrence of large correlations. For independent Gaussian random walks, \citet{Ernst2023} obtained an exact second-moment formula and rates of convergence to the corresponding Wiener-process correlation.

Measures that detect nonlinear dependence raise a related question: does nonsense correlation persist when Pearson correlation is replaced by a coefficient that characterizes independence? Distance correlation is a natural example. For random vectors with finite first moments, its population value is zero if and only if the vectors are independent \citep{Szekely2007}. This characterization extends beyond linear association, but its inferential interpretation still depends on the sampling scheme. Under bivariate normality, distance correlation and several rank-based measures, including Hoeffding's $D$, Blum--Kiefer--Rosenblatt's $R$ and Bergsma--Dassios--Yanagimoto's $\tau^*$, are functions of Pearson correlation \citep{Szekely2007,Zhou2024}. Such population relationships do not settle the random-walk problem.

The difficulty is that process independence does not imply independent sampling across time. At each fixed time, two independent random walks have independent levels and hence zero population distance correlation. A sample formed from one pair of trajectories, however, consists of serially dependent observations: successive levels share accumulated increments, and their distributions change with time. Even with Gaussian increments, this sample is not a collection of independent observations from one bivariate normal distribution. Its asymptotic behaviour therefore requires a time-series analysis. For stationary time series, distance correlation was studied by \citet{Davis2018}. Under stationarity and ergodicity, \citet{Kroll2022} established strong consistency of empirical distance covariance, with null distribution theory under additional dependence conditions. \citet{Betken2024} developed a bootstrap test for independence of strictly stationary time series. These stationary results do not cover random walks.

A separate line of work defines distance covariance between entire stochastic processes \citep{Matsui2017}. Such a coefficient concerns the joint law of the processes as random elements. Here the statistic is computed from time-indexed scalar observations along a single pair of trajectories. The relevant distribution is obtained by sampling time along those realized paths. The resulting coefficient consequently describes association within a path pair, rather than dependence between processes viewed as random elements.

Numerical evidence for nonsense distance correlation appears in \citet[\S\,21.1, Example~21.1]{SzekelyRizzo2023}. Their simulations of 1000 pairs of independent symmetric random walks, each of length 500, produced substantial sample distance correlations. They attributed this behaviour to serial dependence and nonstationarity. We develop a limiting-distribution theory for that phenomenon and establish its consequences for ordinary permutation calibration.

Our analysis assumes centred random walks with independent, identically distributed increments of finite, positive variance. It starts from an exact representation of the sample coefficient as the distance correlation of an empirical occupation measure. After rescaling, the walks converge to independent Brownian paths, and their empirical occupation measure has a random path-based limit. Conditional on the two limiting paths, this measure is obtained by sampling a common time uniformly on $[0,1]$ \citep[p.\,110]{MortersPeres2010}. Independence of the processes does not force this conditional measure to factor into its marginals. This distinction allows persistent sample association despite independence of the original processes.

We establish three results. First, sample distance correlation converges in distribution to a universal Brownian functional that lies strictly between zero and one almost surely and has full support on $[0,1]$. In particular, the probability of exceeding any fixed threshold in $(0,1)$ remains bounded away from zero. Second, we calculate the expected squared distance covariance of the limiting occupation measure. Third, an exact identity for the conditional mean of the permuted statistic yields rejection of process independence with probability tending to one at every fixed nominal level. Together, these results characterize the persistence of nonsense distance correlation and quantify its inferential consequences.

Section~\ref{sec:limit} develops the occupation-measure representation and its limiting properties. Section~\ref{sec:permutation} establishes the failure of ordinary permutation calibration, and Section~\ref{sec:simulation} describes the simulation design and interprets the results for random walks and innovations. Section~\ref{sec:discussion} discusses the implications and open questions. Proofs are collected in Appendix~\ref{app:proofs}; additional simulation details are given in Appendix~\ref{app:simulation}.

\section{The occupation-measure limit}
\label{sec:limit}

\subsection{Population, sample and path representations}
\label{subsec:representations}

We first express population, sample and path distance correlations through the same centring operation. For real random variables $X,Y$ satisfying $\E(|X|)+\E(|Y|)<\infty$, let $(X',Y')$ be an independent copy of $(X,Y)$. Define the population-centred $X$-distance \citep{Szekely2007,SzekelyRizzo2009}
$$
 A(x,x')=|x-x'|-\E(|x-X|)-\E(|x'-X|)+\E(|X-X'|),
$$
and define $B(y,y')$ analogously for $Y$. The squared population distance covariance and correlation admit the representation \citep[\S\,2]{Lyons2013}
$$
 \begin{aligned}
 V^2(X,Y)&=\E\{A(X,X')B(Y,Y')\},\\
 R^2(X,Y)&=\frac{V^2(X,Y)}{\{V^2(X,X)V^2(Y,Y)\}^{1/2}},
 \end{aligned}
$$
with $R(X,Y)$ the nonnegative square root and $R(X,Y)=0$ when the denominator vanishes. Here the expectation is over two independent copies of the joint pair; independence between $X$ and $Y$ is not assumed. Distance correlation lies in $[0,1]$ and equals zero if and only if $X$ and $Y$ are independent. It is invariant under separate translations and nonzero rescalings.

The sample version replaces the expectations in $A$ and $B$ by empirical averages. For real observation pairs $(x_i,y_i)$, let
\[
 a_{ij}=|x_i-x_j|,\qquad
 A_{ij}=a_{ij}-\overline{a}_{i\cdot}-\overline{a}_{\cdot j}+\overline{a}_{\cdot\cdot}
 \quad(i,j=1,\ldots,n),
\]
where the row, column and overall means are
\[
 \overline{a}_{i\cdot}=\frac1n\sum_{j=1}^n a_{ij},\qquad
 \overline{a}_{\cdot j}=\frac1n\sum_{i=1}^n a_{ij},\qquad
 \overline{a}_{\cdot\cdot}=\frac1{n^2}\sum_{i,j=1}^n a_{ij}.
\]
Define the distances $b_{ij}$, their row, column and overall means, and the centred distances $B_{ij}$ analogously for the $y$ observations. For the centred distance matrices $A=(A_{ij})$ and $B=(B_{ij})$, define their normalized inner product by
\[
 \langle A,B\rangle_n=\frac1{n^2}\sum_{i,j=1}^n A_{ij}B_{ij}.
\]
The usual sample distance correlation, based on empirical double centring, satisfies
\begin{equation}
 R_n^2=\frac{\langle A,B\rangle_n}
 {\{\langle A,A\rangle_n\langle B,B\rangle_n\}^{1/2}}.
 \label{eq:sample}
\end{equation}
Set $R_n=0$ when the denominator vanishes, and otherwise take the nonnegative square root. All sums include the diagonal terms. This coefficient is the distance correlation of the empirical distribution assigning mass $1/n$ to each pair, even for dependent observations.

For bounded measurable paths $f,g:[0,1]\to\mathbb R$, define their joint occupation measure on Borel sets $H\subseteq\mathbb R^2$ by
\[
 \mu_{f,g}(H)=\int_0^1\ind{(f(t),g(t))\in H}\dd t.
\]
Equivalently, for $U$ uniform on $[0,1]$, $\mu_{f,g}$ is the law of $(f(U),g(U))$: it assigns to each region the fraction of time the paired path spends there \citep[p.\,110]{MortersPeres2010}. Both coordinates use the same time $U$, preserving the observed pairing. Sampling them at independent uniform times would instead produce the product of the two marginal occupation measures.

For a path $f$, let $D_f(s,t)$ denote its pairwise distance centred with respect to uniform time:
$$
 \begin{aligned}
 D_f(s,t)&=|f(s)-f(t)|-\int_0^1|f(s)-f(u)|\dd u\\
 &\quad-\int_0^1|f(t)-f(u)|\dd u+\int_0^1\!\int_0^1|f(u)-f(v)|\dd u\dd v.
 \end{aligned}
$$
Define $D_g$ by the same formula with $g$ in place of $f$. These are the path analogues of the centred sample distances $A_{ij}$ and $B_{ij}$. The squared distance covariance $K(f,g)=V^2\{f(U),g(U)\}$ and distance correlation $\mathcal R(f,g)=R\{f(U),g(U)\}$ admit the integral forms
\begin{equation}
 \begin{aligned}
 K(f,g)&=\int_0^1\!\int_0^1 D_f(s,t)D_g(s,t)\dd s\dd t,\\
 \mathcal R^2(f,g)&=\frac{K(f,g)}{\{K(f,f)K(g,g)\}^{1/2}}.
 \end{aligned}
 \label{eq:path}
\end{equation}
The squared distance variances in the denominator are
\[
 K(f,f)=\int_0^1\!\int_0^1D_f(s,t)^2\dd s\dd t,
 \qquad
 K(g,g)=\int_0^1\!\int_0^1D_g(s,t)^2\dd s\dd t.
\]
The two integration variables $s,t$ correspond to two independent uniform times used to draw independent pairs from $\mu_{f,g}$. We set $\mathcal R(f,g)=0$ when the denominator vanishes. Then $0\leq\mathcal R(f,g)\leq1$, and $K(f,f)>0$ precisely when $f$ is not constant almost everywhere. In particular, the denominator is positive for two continuous nonconstant paths.

For random paths, the definitions are applied to each realization, with $U$ independent of the processes. The occupation measure is then random, and $K$ and $\mathcal R$ average over time conditional on the paths. Process independence concerns the distribution across realizations; the occupation coefficient concerns the distribution across time within a realization. Keeping these two sources of randomness distinct is essential to interpreting the limit.

For step paths taking values $x_i,y_i$ on $((i-1)/n,i/n]$, the occupation measure $\mu_{f,g}$ assigns mass $1/n$ to each observed pair. The integrals in \eqref{eq:path} then equal the corresponding sample averages, so $\mathcal R(f,g)=R_n$ exactly. Under the assumptions below, weak convergence of the rescaled walks to Brownian paths and continuity of $\mathcal R$ at pairs of nonconstant continuous paths yield the limiting distribution of $R_n$. Appendix~\ref{app:proofs} establishes the required continuity and limiting properties.

\subsection{A nondegenerate universal limit}
\label{subsec:limit}

Let $(\xi_j)$ and $(\eta_j)$ be mutually independent sequences, each consisting of independent and identically distributed real random variables. Assume
\[
 \E(\xi_1)=\E(\eta_1)=0,\qquad
 0<\sigma_\xi^2=\E(\xi_1^2)<\infty,\qquad
 0<\sigma_\eta^2=\E(\eta_1^2)<\infty.
\]
Define the two random walks by
\[
 X_0=Y_0=0,\qquad
 X_k=\sum_{j=1}^k\xi_j,\qquad
 Y_k=\sum_{j=1}^k\eta_j,\qquad k=1,2,\ldots,
\]
and compute $R_n$ from the first $n$ paired levels $(X_k,Y_k)$. All asymptotic statements concern the sample size $n\to\infty$. On the time scale $k/n$, divide the walks by $\sigma_\xi n^{1/2}$ and $\sigma_\eta n^{1/2}$, respectively. These separate rescalings leave distance correlation unchanged. Donsker's theorem gives convergence in distribution of their interpolated paths to two independent standard Brownian motions on $[0,1]$ \citep[Theorem~5.22]{MortersPeres2010}.

Write $W_1$ and $W_2$ for these independent standard Brownian motions on $[0,1]$. Applying the path functional in \eqref{eq:path} defines
\[
 R_\infty=\mathcal R(W_1,W_2).
\]
Thus $R_\infty$ varies across path pairs through their joint occupation measure.

\begin{theorem}
\label{thm:limit}
Under the preceding innovation assumptions, as $n\to\infty$,
\[
 R_n\dto R_\infty.
\]
The limiting random variable satisfies
\[
 \pr(0<R_\infty<1)=1,\qquad
 \supp\{\law(R_\infty)\}=[0,1].
\]
Consequently, for every fixed $c\in(0,1)$ and every fixed positive integer $q\in\{1,2,\ldots\}$,
\[
 \begin{aligned}
 \liminf_{n\to\infty}\pr(R_n>c)&\geq\pr(R_\infty>c)>0,\\
 \lim_{n\to\infty}\E(R_n^q)&=\E(R_\infty^q)>0.
 \end{aligned}
\]
\end{theorem}

The theorem describes how the distribution of $R_n$ across repeated samples changes as $n$ increases. It does not assert almost sure convergence along a single successively extended pair of trajectories. Full support on $[0,1]$ means that every nonempty open subinterval of $(0,1)$ has positive probability, even though neither endpoint is attained with positive probability. The integer $q$ is held fixed as $n\to\infty$; the choices $q=1$ and $q=2$ give convergence of the first and second moments. The nondegenerate limit rules out an explanation based solely on small-sample bias. The functional central limit theorem and scale invariance give universality.

The positive limit concerns the occupation measure obtained by fixing the paths and sampling a common time. Zero distance correlation under this measure would imply independence of the paired values and hence zero time covariance. However, conditional on $W_1$, the time covariance with $W_2$ is a nondegenerate Gaussian variable. It is therefore nonzero almost surely. The Appendix supplies this argument, the continuity details, and the support construction.

We also obtain an explicit expectation for the squared distance covariance in \eqref{eq:path}.

\begin{proposition}
\label{prop:moment}
Let $W_1$ and $W_2$ be the independent standard Brownian motions on $[0,1]$ defined above, and let $K$ be the squared distance covariance functional in \eqref{eq:path}. Then
\[
 \E\{K(W_1,W_2)\}=\frac{26}{75\pi}-\frac{1}{12}=0.0270140939\ldots.
\]
\end{proposition}

Proposition~\ref{prop:moment} averages the pathwise squared distance covariance over independent pairs of Brownian paths. It concerns the unnormalized numerator in \eqref{eq:path}; it does not give $\E(R_\infty^2)$, because the denominator is random and cannot be replaced by its expectation.

\section{Failure of ordinary permutation calibration}
\label{sec:permutation}

The occupation limit also has consequences for permutation calibration. Consider the squared sample distance covariance
\[
 T_n=\langle A,B\rangle_n.
\]
The overall distance means $\overline{a}_{\cdot\cdot}$ and $\overline{b}_{\cdot\cdot}$ retain their definitions from Section~\ref{subsec:representations}. Let $\pi$ be uniformly distributed over the permutations of $1,\ldots,n$, independently of the data. Permuting the time labels of the second series gives the statistic
\[
 T_n^{\pi}=\frac1{n^2}\sum_{i,j=1}^n A_{ij}B_{\pi(i)\pi(j)}.
\]
Ordinary permutation calibration uses the exact conditional tail probability
\begin{equation}
 p_n=\pr_\pi(T_n^{\pi}\geq T_n\mid\mathcal D_n),
 \label{eq:pvalue}
\end{equation}
where $\mathcal D_n$ denotes the observed sample. If $T_n=0$, this probability is one. Since the denominator in \eqref{eq:sample} is unchanged by permutation, \eqref{eq:pvalue} also calibrates the corresponding upper-tail distance-correlation test whenever the denominator is nonzero.

\begin{theorem}
\label{thm:permutation}
For any sample of size $n\geq2$,
\begin{equation}
 \E_\pi(T_n^{\pi}\mid\mathcal D_n)=\frac{\overline{a}_{\cdot\cdot} \overline{b}_{\cdot\cdot}}{n-1}.
 \label{eq:permmean}
\end{equation}
On the event $T_n>0$, the exact permutation probability satisfies
\[
 p_n\leq\frac{\overline{a}_{\cdot\cdot} \overline{b}_{\cdot\cdot}}{(n-1)T_n}.
\]
Under the assumptions of Theorem~\ref{thm:limit}, $p_n\prto0$ and $\pr(p_n\leq\alpha)\to1$ for every fixed $\alpha\in(0,1)$.
\end{theorem}

Theorem~\ref{thm:permutation} quantifies a failure of calibration under process independence. Ordinary permutation changes the temporal organization of the observations, which are not exchangeable under the random-walk null. The resulting reference distribution therefore does not represent that null. The conclusion applies to unrestricted permutation of time labels; it makes no claim about resampling procedures designed to preserve temporal dependence.

\section{Numerical results}
\label{sec:simulation}

\subsection{Simulation design and reported summaries}

The experiments address three questions: whether distance correlation persists for independent random walks, how its distribution changes with sample size, and whether ordinary permutation calibration maintains its nominal rejection rate. We use Gaussian innovations to illustrate the finite-variance theory, rather than to impose an additional assumption on Theorem~\ref{thm:limit}. 

For each sample size $n\in\{50,100,250,500,1000,2000,5000\}$, we generated $M=2000$ independent pairs of innovation sequences. Within each replication, all $2n$ innovations were independent standard normal variables. Their cumulative sums form the random walks defined in Section~\ref{subsec:limit}. We computed \eqref{eq:sample} first from the $n$ pairs $(X_k,Y_k)$ and then from the original pairs $(\xi_k,\eta_k)$. The two inputs are labelled \emph{Random walks} and \emph{Innovations}, respectively. They share the same innovations within a replication, so the comparison isolates the effect of accumulation. Fresh sequences were generated for each replication and sample size; larger samples were not extensions of the smaller simulated paths. 

Each replication produces one distance correlation for each input. Table~\ref{tab:simulation} summarizes these coefficients across the $M$ replications. Its mean columns report the average $R_n$, not the square root of the average $R_n^2$. The quantile columns give the $0.050$ and $0.950$ empirical quantiles of the random-walk coefficients, so they delimit the central 90\% of their simulated distribution. They are neither quantiles of the observations $X_k,Y_k$ nor confidence limits for the mean. The last column is the fraction of replications with random-walk correlation above $0.300$. These summaries separate the typical size of the coefficient from its variation across independent pairs of trajectories.

\begin{table}
\centering
\caption{Sample distance correlation for independent Gaussian sequences}
\small
\begin{tabular}{rrrrrr}
 & \multicolumn{2}{c}{Mean} & \multicolumn{2}{c}{Random-walk $R_n$ quantile} & \\
 $n$ & Random walks & Innovations & 0.050 & 0.950 & $\pr(R_n>0.300)$\\
50 & 0.510 & 0.244 & 0.247 & 0.854 & 0.860 \\
100 & 0.504 & 0.173 & 0.233 & 0.843 & 0.834 \\
250 & 0.508 & 0.110 & 0.218 & 0.857 & 0.822 \\
500 & 0.512 & 0.078 & 0.217 & 0.857 & 0.833 \\
1000 & 0.499 & 0.055 & 0.220 & 0.853 & 0.824 \\
2000 & 0.503 & 0.039 & 0.218 & 0.848 & 0.829 \\
5000 & 0.500 & 0.025 & 0.215 & 0.844 & 0.822 \\
\end{tabular}
\label{tab:simulation}
\par\smallskip\begin{minipage}{\textwidth}\footnotesize
Each row uses 2000 replications. Quantiles and exceedance probabilities concern only the random-walk input. The quantile range describes variation between replications; standard errors of the means are reported separately in Table~\ref{tab:precision}.
\end{minipage}
\end{table}

\subsection{Persistent correlations and interpretation of the figure}

The mean distance correlation for random walks remains near $0.500$, whereas the mean for innovations decreases as $n$ grows. At $n=5000$, Table~\ref{tab:simulation} gives means of $0.500$ and $0.025$, respectively, and random-walk quantiles of $0.215$ and $0.844$. Thus even a long pair of independent random walks can produce a substantial coefficient. The estimated exceedance probability is $0.822$: 1644 of the 2000 random-walk pairs have correlation above $0.300$. The value near $0.500$ is an empirical finding, not an analytic value established for the mean of $R_\infty$.

The contrast has a direct sampling interpretation. The innovation pairs are independent observations from a fixed product distribution, so their empirical distribution approaches that independent population as the sample grows. Random-walk observations instead reuse earlier increments: one shock affects every subsequent value of its trajectory. Increasing $n$ therefore adds temporally related observations rather than independent copies of a fixed pair of variables. Long excursions can create persistent patterns of association within a realized pair of paths, even though the processes generating those paths are independent.

The pairwise distances make this mechanism more explicit. For the standard Gaussian random walks and distinct time indices $i,j$,
\[
 \E(|X_i-X_j|)=\E(|Y_i-Y_j|)=\left(\frac{2}{\pi}\right)^{1/2}|i-j|^{1/2}.
\]
Both distance matrices therefore have a systematic dependence on the separation of the observation times. This common time structure can contribute to their association without any dependence between the two processes. The occupation-measure limit gives the rigorous explanation: after rescaling time and amplitude, the statistic approaches a random functional of two Brownian paths, rather than the distance correlation of an independent pair sampled at a fixed time. Scale invariance removes the amplitude normalization, while the realized temporal patterns remain. This explains persistence without specifying the limiting mean.

\begin{figure}
\centering
\includegraphics[width=\textwidth]{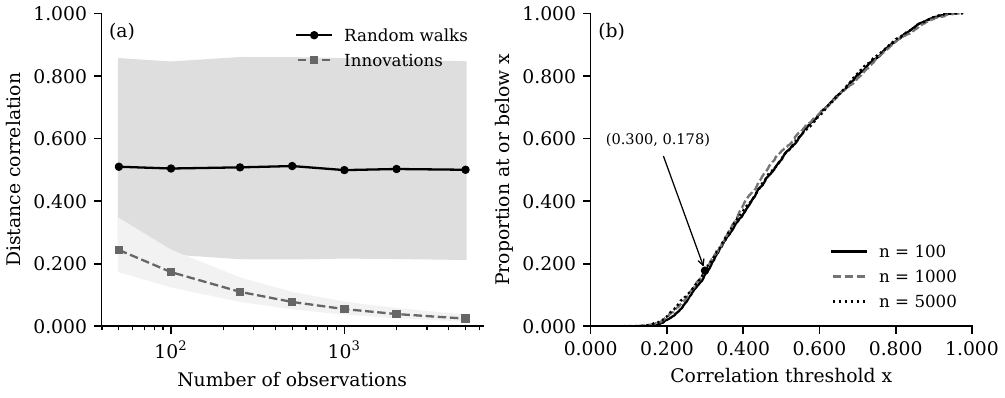}
\caption{Distance correlation for independent Gaussian inputs. Panel (a): means for random walks (solid line, circles) and innovations (dashed line, squares), with their $0.050$--$0.950$ quantile bands (dark and light shading). Panel (b): empirical cumulative distribution functions of random-walk correlations for $n=100$ (solid), $n=1000$ (dashed) and $n=5000$ (dotted). The marked point gives the fraction of the $n=5000$ coefficients not exceeding $0.300$.}
\label{fig:simulation}
\par\smallskip\noindent Alt text: Random-walk correlations remain dispersed as sample size increases. Their cumulative distribution curves nearly overlap; the marked point is $(0.300,0.178)$.
\end{figure}

Figure~\ref{fig:simulation}(a) displays the same mean comparison over all sample sizes. Its shaded bands connect the $0.050$ and $0.950$ quantiles for each input. The random-walk band stays broad, whereas the innovation band moves towards zero. A broad quantile band is compatible with a precisely estimated mean: at $n=5000$, the random-walk coefficients have standard deviation $0.196$, but the Monte Carlo standard error of their mean is only $4.394\times10^{-3}$. The former measures variation between path pairs; the latter measures simulation uncertainty in the reported average.

Panel (b) compares entire distributions rather than means. If $R_n^{(r)}$ denotes the random-walk coefficient from replication $r$, the plotted empirical cumulative distribution function is
\[
 \widehat F_n(x)=\frac1M\sum_{r=1}^M I(R_n^{(r)}\leq x),\qquad M=2000.
\]
The three curves are close over a broad range of thresholds, indicating similar distributions at the displayed sample sizes. If the coefficients $R_n^{(r)}$ were concentrating at zero, these curves would approach  one at every fixed positive threshold. Their persistent spread is consistent with the nondegenerate limit in Theorem~\ref{thm:limit}; visual agreement alone does not establish a convergence rate.

\subsection{Permutation calibration and Monte Carlo uncertainty}

We assessed calibration in a separate experiment with 300 new replications at each of $n=100$, 500 and 2000. For each input, we held the first series fixed and independently permuted the time labels of the second series 199 times. Each permutation produces a new distance correlation. If $k$ of these coefficients are at least as large as the observed coefficient, the Monte Carlo permutation probability is $(k+1)/200$. We rejected independence when this probability was at most $0.050$. The plus-one adjustment includes the observed arrangement and gives a minimum possible probability of $0.005$; Theorem~\ref{thm:permutation} concerns the exact probability over all permutations.

For random walks, the rejection proportions were $0.940$, $1.000$ and $1.000$ at the three respective sample sizes. For innovations, they were $0.050$, $0.047$ and $0.050$. Table~\ref{tab:permutation} reports the counts and confidence intervals. Since the two processes are independent in every replication, these are false-rejection frequencies, not power against a dependent alternative. The innovation results are close to the nominal level. Permuting random-walk observations, however, breaks their temporal organization and produces a reference distribution that does not represent independent random walks. The resulting mismatch explains the large rejection frequencies.

A reported proportion is itself estimated from a finite number of simulations. We quantify this uncertainty using two-sided 95\% Wilson score intervals for a binomial probability \citep{Wilson1927}. Here a trial is one independent simulation, and a success means either rejection or exceedance of a stated correlation threshold. For example, the exceedance estimate $0.822$ has interval $(0.805,0.838)$ based on 2000 replications. This interval quantifies uncertainty in the exceedance probability, whereas Table~\ref{tab:simulation} describes variation among coefficients.

For 300 rejections in 300 trials, the Wilson interval is $(0.987,1.000)$ to three decimal places. The value $0.987$ is the lower endpoint of the two-sided 95\% interval; a one-sided 95\% confidence bound would use a different normal quantile. Observing only rejections does not establish that the underlying rejection probability is exactly one; the interval records uncertainty from the finite simulation count. The 95\% level describes the approximate long-run coverage of the interval procedure under repeated simulation experiments. The formula and its all-rejections specialization are given in Appendix~\ref{app:precision}.

\subsection{Covariance check and implementation}

As a further check on Proposition~\ref{prop:moment}, we divided both random-walk inputs by $n^{1/2}$ and computed their squared distance covariance. This rescaling leaves distance correlation unchanged but divides squared distance covariance by $n$, giving a quantity with a finite Brownian limit. At $n=5000$, its simulated mean is $26.147\times10^{-3}$ with Monte Carlo standard error $0.804\times10^{-3}$, compared with the theoretical limiting expectation $27.014\times10^{-3}$. Their difference is about $1.078$ Monte Carlo standard errors. This agreement is consistent with the limit calculation, while allowing for both Monte Carlo variability and finite-sample approximation error.

\section{Discussion}
\label{sec:discussion}

The results distinguish two inferential targets: dependence between stochastic processes and association under the time distribution of a realized path pair. This distinction matters when nonstationary trajectories are analysed through observation pairs. In the present model, first differencing recovers independent innovations, but changes the target from levels to increments. A test designed for increments consequently answers a different question unless the model connects that question to process independence.

Valid calibration for nonstationary observations must account for their temporal dependence under the null. The bootstrap theory of \citet{Betken2024} shows how this can be achieved under explicit stationary dependence assumptions; it does not by itself justify resampling random-walk levels. Constructing a calibrated procedure for the nonstationary null is therefore a separate problem from identifying the limiting occupation coefficient.

Several distributional questions also remain. Theorem~\ref{thm:limit} identifies the limiting functional, its support and the convergence of its moments, but does not provide a density, explicit correlation moments or finite-sample approximation rates. Proposition~\ref{prop:moment} evaluates the expected numerator before normalization; the random distance variances make the corresponding correlation moments more difficult. Extensions to infinite-variance innovations, long-memory partial sums and other nonlinear dependence measures require their own functional limits and continuity arguments.

\appendix
\section{Proofs}
\label{app:proofs}

Throughout the proofs, write
\begin{align*}
 a_f(s,t)&=|f(s)-f(t)|, & r_f(s)&=\int_0^1a_f(s,u)\dd u, & m_f&=\int_0^1r_f(u)\dd u,\\
 D_f(s,t)&=a_f(s,t)-r_f(s)-r_f(t)+m_f.
\end{align*}
Here $r_f$ and $m_f$ are the averages over one and two time arguments, respectively. Define $a_g,r_g,m_g$ analogously, so that $D_g(s,t)=a_g(s,t)-r_g(s)-r_g(t)+m_g$.

\begin{proof}[Proof of Theorem~\ref{thm:limit}]
We first establish continuity of the functionals. If $f_l\to f$ uniformly, then
\[
 \|a_{f_l}-a_f\|_\infty\leq2\|f_l-f\|_\infty,
 \qquad \|D_{f_l}-D_f\|_\infty\leq8\|f_l-f\|_\infty.
\]
Thus $K$ is continuous under uniform convergence on bounded paths. For a continuous $f$, the identity $K(f,f)=0$ implies $D_f(s,t)=0$ for every $s,t$. Taking $s=t$ gives $r_f(s)=m_f/2$; integration gives $m_f=0$, and hence $f$ is constant. Therefore $\mathcal R$ is continuous at every pair of continuous nonconstant paths.

Let $\widetilde X_n$ be the polygonal interpolation through $X_k/(\sigma_\xi n^{1/2})$ at $k/n$, including $X_0=0$, and define $\widetilde Y_n$ analogously. Donsker's theorem gives joint weak convergence in $C_0[0,1]^2$ to $(W_1,W_2)$; see \citet[Theorem~5.22]{MortersPeres2010}. On the interval $((k-1)/n,k/n]$, let $\widehat X_n$ equal $X_k/(\sigma_\xi n^{1/2})$ and define $\widehat Y_n$ analogously, with $\widehat X_n(0)=\widehat Y_n(0)=0$. Then
\[
 \|\widehat X_n-\widetilde X_n\|_\infty
 \leq\max_{1\leq k\leq n}|\xi_k|/(\sigma_\xi n^{1/2})=o_{\pr}(1),
\]
and the corresponding statement holds for $Y$. The last assertion follows from finite second moments, since $n\pr(|\xi_1|>\epsilon n^{1/2})\to0$ for every $\epsilon>0$. Integrating the step functions gives exactly the discrete averages in \eqref{eq:sample}. Scale invariance, the continuity just established, and the fact that Brownian paths are nonconstant almost surely prove $R_n\dto\mathcal R(W_1,W_2)$.

To prove strict positivity, first observe that if $\mathcal R(W_1,W_2)=0$, the characterization of independence by distance correlation implies zero covariance under $\mu_{W_1,W_2}$. Consequently,
\[
 Z=\int_0^1\{W_1(t)-\overline{W}_1\}W_2(t)\dd t=0,
 \qquad \overline{W}_1=\int_0^1 W_1(t)\dd t.
\]
Conditionally on $W_1$, stochastic Fubini gives
\[
 \begin{aligned}
 F_{W_1}(u)&=\int_u^1\{W_1(t)-\overline{W}_1\}\dd t,\\
 Z&=\int_0^1F_{W_1}(u)\dd W_2(u),\qquad
 \var(Z\mid W_1)=\int_0^1F_{W_1}(u)^2\dd u.
 \end{aligned}
\]
The stochastic integral is conditionally Gaussian. Its variance is positive almost surely: if it vanished, continuity would give $F_{W_1}\equiv0$, whose derivative forces $W_1$ to be constant. It follows that $\pr(Z=0\mid W_1)=0$ almost surely, hence $\pr\{\mathcal R(W_1,W_2)=0\}=0$.

For the strict upper bound, equality in Cauchy--Schwarz implies that $\mathcal R(f,g)=1$ requires $D_g=cD_f$ for some $c>0$, initially almost everywhere and then everywhere by continuity. Evaluation on the diagonal and integration yield $m_g=cm_f$ and $r_g=cr_f$. Thus $|g(s)-g(t)|=c|f(s)-f(t)|$ for every $s,t$. The map from the interval $f([0,1])$ to $g([0,1])$ induced by the paths is a scaled isometry, so $g=a+bf$ with $b=\pm c$. Since the Brownian paths start at zero, this would require $W_2=bW_1$. In particular,
\[
 W_1(1/2)W_2(1)-W_1(1)W_2(1/2)=0.
\]
The vector of four evaluations has a nondegenerate Gaussian density, and the displayed nonzero polynomial has a zero set of Lebesgue measure zero. Hence $\pr\{\mathcal R(W_1,W_2)=1\}=0$.

For the support assertion, Wiener measure has full support on $C_0[0,1]$; see \citet[Exercise~1.8]{MortersPeres2010}. This property follows by approximating a continuous path by a polygonal one and using the positive probabilities of suitable Gaussian increments and Brownian-bridge tubes. Independence gives full support for the path pair. The choice $f(t)=g(t)=t$ yields $\mathcal R(f,g)=1$.

To obtain arbitrarily small correlations, let $h$ be the continuous triangular wave of period one, with $h(v)=2v$ for $0\leq v\leq1/2$ and $h(v)=2(1-v)$ for $1/2\leq v\leq1$. Take $f(t)=t$ and $g_m(t)=h(mt)$, where $m$ is a positive integer. Both occupation marginals are uniform on $[0,1]$. For any continuous $\phi$ on $[0,1]^2$,
\begin{align*}
 \int_0^1\phi\{t,g_m(t)\}\dd t
 &=m^{-1}\sum_{j=0}^{m-1}\int_0^1\phi\{(j+v)/m,h(v)\}\dd v\\
 &\longrightarrow\int_0^1\!\int_0^1\phi(x,y)\dd x\dd y.
\end{align*}
The convergence follows from a uniform Riemann-sum approximation and the uniform distribution of $h(U)$. Since all measures are supported on a common compact set, the distance-moment representation is continuous under this weak convergence. Thus $K(f,g_m)\to0$, whereas $K(f,f)$ and $K(g_m,g_m)$ equal the same positive constant. Therefore $\mathcal R(f,g_m)\to0$.

For any $0<r<1$, choose $m$ such that $\mathcal R(f,g_m)<r$. The interpolation $g_\lambda=(1-\lambda)f+\lambda g_m$, $0\leq\lambda\leq1$, is nonconstant throughout. Indeed, $g_\lambda(0)=0$ and $g_\lambda(1)=1-\lambda$ when $\lambda<1$, while $g_1=g_m$ is nonconstant. Continuity and the intermediate value theorem give a path pair with correlation $r$. Full support of Wiener measure and continuity of $\mathcal R$ imply that every open interval containing $r$ receives positive probability. The two endpoint constructions establish support $[0,1]$.

The tail assertion follows from the Portmanteau theorem. The moment assertion follows from weak convergence and $0\leq R_n\leq1$.
\end{proof}

\begin{proof}[Proof of Proposition~\ref{prop:moment}]
Independence of the Brownian paths gives
\[
 \E\{D_{W_1}(s,t)D_{W_2}(s,t)\}=
 \E\{D_{W_1}(s,t)\}\E\{D_{W_2}(s,t)\}.
\]
Since $\E(|W_1(s)-W_1(t)|)=(2/\pi)^{1/2}|s-t|^{1/2}$, let $d(s,t)=|s-t|^{1/2}$ and set
\[
 r(s)=\int_0^1d(s,t)\dd t=\frac23\{s^{3/2}+(1-s)^{3/2}\},
 \qquad m=\int_0^1r(s)\dd s=\frac8{15}.
\]
The required integrals are
\[
 \int_0^1\!\int_0^1d(s,t)^2\dd s\dd t=\frac13,
 \qquad \int_0^1r(s)^2\dd s=\frac29+\frac\pi{48}.
\]
Expansion of the squared centred kernel therefore gives
\[
 \E\{K(W_1,W_2)\}
 =\frac2\pi\left\{\frac13-2\left(\frac29+\frac\pi{48}\right)
 +\left(\frac8{15}\right)^2\right\}
 =\frac{26}{75\pi}-\frac1{12}.
\]
Finite second moments of the Brownian suprema justify the interchanges of expectation and integration.
\end{proof}

\begin{proof}[Proof of Theorem~\ref{thm:permutation}]
Double centring gives zero row sums, diagonal sum $-n\overline{a}_{\cdot\cdot}$ and off-diagonal sum $n\overline{a}_{\cdot\cdot}$ for $A$, with analogous identities for $B$. For a uniform random permutation, the expected diagonal entry of the permuted $B$ is $-\overline{b}_{\cdot\cdot}$ and the expected off-diagonal entry is $\overline{b}_{\cdot\cdot}/(n-1)$. Summing separately over diagonal and off-diagonal indices gives \eqref{eq:permmean}. Conditional Markov's inequality applied to the nonnegative statistic $T_n^{\pi}$ gives the bound for $T_n>0$.

Normalize the two samples by their respective innovation standard deviations times $n^{1/2}$. The ratio $\overline{a}_{\cdot\cdot} \overline{b}_{\cdot\cdot}/T_n$ is unchanged. The continuity argument used for Theorem~\ref{thm:limit} yields joint weak convergence of the normalized average distances and test statistic to $m_{W_1}$, $m_{W_2}$ and $K(W_1,W_2)$. Strict positivity of $K(W_1,W_2)$ implies that $\overline{a}_{\cdot\cdot} \overline{b}_{\cdot\cdot}/T_n$ is bounded in probability on $T_n>0$, and that $\pr(T_n=0)\to0$. The permutation bound gives $p_n\prto0$ and hence the fixed-level rejection assertion.
\end{proof}

\section{Additional simulation details}
\label{app:simulation}
\subsection{Simulation design}

All results use the original nonnegative sample distance correlation based on doubly centred Euclidean distance matrices. Reported correlations are $R_n$, not $R_n^2$. For each sample size, two arrays of independent standard normal innovations were generated. The input labelled \emph{Random walks} uses their cumulative sums; \emph{Innovations} uses the increments themselves. Each sample size used 2000 replications generated from independent innovation arrays.

The master seed was 20261110. NumPy's \texttt{SeedSequence} generated three child streams for implementation checks, main simulations and permutation experiments. Each used \texttt{default\_rng} with its default PCG64 generator. The recorded environment was Python 3.12.14 with NumPy 2.5.4, SciPy 1.18.1, Matplotlib 3.11.2, \texttt{dcor} 0.7 and Numba 0.68.1. All data were simulated; no external dataset was used.

The statistic was computed with version 0.7 of the \texttt{dcor} package and its mergesort method. For one-dimensional inputs this avoids storing an entire pair of distance matrices. At sample sizes 7, 20, 53 and 101, both innovation and walk examples were also evaluated by explicit matrix double centring. The largest absolute discrepancy was $3.886\times10^{-15}$. Additional checks confirmed translation and scale invariance, and unit correlation for nonconstant affine relationships.

\subsection{Monte Carlo precision}
\label{app:precision}

For simulated coefficients $R_n^{(1)},\ldots,R_n^{(M)}$, the standard error of their average was the sample standard deviation divided by $M^{1/2}$. Table~\ref{tab:precision} reports standard errors for the mean correlations and for the squared distance covariance after rescaling. The quantity in the final two columns is the squared sample distance covariance after both Gaussian-walk inputs have been divided by $n^{1/2}$. Its Brownian limiting expectation is $26/(75\pi)-1/12$, approximately $27.014\times10^{-3}$.

\begin{table}
\centering
\caption{Monte Carlo precision (entries in units of $10^{-3}$)}
\small
\begin{tabular}{rrrrr}
 & \multicolumn{2}{c}{Standard error of mean correlation} & \multicolumn{2}{c}{Scaled squared distance covariance}\\
 $n$ & Random walks & Innovations & Mean & Standard error\\
50 & 4.237 & 1.174 & 27.481 & 0.895 \\
100 & 4.328 & 0.791 & 25.944 & 0.777 \\
250 & 4.509 & 0.503 & 28.529 & 0.906 \\
500 & 4.459 & 0.352 & 28.387 & 0.877 \\
1000 & 4.462 & 0.249 & 26.589 & 0.846 \\
2000 & 4.424 & 0.183 & 26.482 & 0.823 \\
5000 & 4.394 & 0.108 & 26.147 & 0.804 \\
\end{tabular}
\label{tab:precision}
\par\smallskip\begin{minipage}{\textwidth}\footnotesize
Each row uses 2000 replications. All four numerical columns after $n$ are expressed in units of $10^{-3}$ and rounded to three decimal places; for example, $4.394$ represents $4.394\times10^{-3}$. The final column describes Monte Carlo error, not finite-sample approximation error. The quantile bands in Section~\ref{sec:simulation} describe the distribution of the statistic, not confidence intervals for its expectation.
\end{minipage}
\end{table}

For an estimated event probability, let $k$ of $M$ independent replications satisfy the event and put $\widehat\theta=k/M$. The two-sided 95\% Wilson interval used in Section~\ref{sec:simulation} has endpoints
\[
 \frac{\widehat\theta+z^2/(2M)
 \ \pm\ z\{\widehat\theta(1-\widehat\theta)/M+z^2/(4M^2)\}^{1/2}}
 {1+z^2/M},\qquad z=\Phi^{-1}(0.975),
\]
where $\Phi$ is the standard normal distribution function. The event is either $R_n>0.300$ or rejection by the permutation test. These intervals describe uncertainty in estimated event probabilities across simulations, not the range of individual correlations. For $k=M$, they simplify to
\[
 \left(\frac{M}{M+z^2},1\right).
\]
For $M=300$, this yields $(0.987,1.000)$ after rounding. The lower endpoint belongs to a two-sided 95\% interval; a one-sided 95\% bound would use a different normal quantile. All calculations use the unrounded value of $z$, and all reported summaries are rounded only after calculation.

\subsection{Permutation experiment}

We deliberately used an independent-sample permutation scheme on the random-walk observations to quantify the resulting failure of calibration. For each sample size, we generated 300 independent datasets and, within each dataset, 199 independent uniform permutations of the second coordinate. Comparing the permuted coefficients with the observed correlation gives
\[
 \widehat p=\frac{1+\sum_{b=1}^{199}
 I(R_n^{\pi_b}\geq R_n)}{200}.
\]
Rejection was defined by $\widehat p\leq0.050$. Its minimum attainable value is $1/200=0.005$. Theorem~\ref{thm:permutation} instead concerns the exact tail probability over all permutations. The corresponding innovation experiment is a valid independent-sample control.

\begin{table}
\centering
\caption{Rejection under ordinary permutation calibration at nominal level 0.050}
\small
\begin{tabular}{rlrrr}
 $n$ & Sample & Rejections & Proportion & Wilson 95\% interval\\
100 & Random walks & 282/300 & 0.940 & (0.907, 0.962) \\
100 & Innovations & 15/300 & 0.050 & (0.031, 0.081) \\
500 & Random walks & 300/300 & 1.000 & (0.987, 1.000) \\
500 & Innovations & 14/300 & 0.047 & (0.028, 0.077) \\
2000 & Random walks & 300/300 & 1.000 & (0.987, 1.000) \\
2000 & Innovations & 15/300 & 0.050 & (0.031, 0.081) \\
\end{tabular}
\label{tab:permutation}
\par\smallskip\begin{minipage}{\textwidth}\footnotesize
Every pair of processes is independent. For random-walk observations, these are rejection proportions from an invalid calibration procedure, not power against process dependence. There are 300 replications and 199 random permutations per replication.
\end{minipage}
\end{table}

\end{document}